\documentclass[journal,twoside,web]{ieeecolor}

\def\logoname{LOGO-generic-web}
\definecolor{subsectioncolor}{rgb}{0,0.541,0.855}
\def\journalname{IEEE Transactions on Automatic Control}

\usepackage{cite}
\usepackage{amsmath,amssymb,amsfonts}
\usepackage{algorithmic}
\usepackage{graphicx}
\usepackage{algorithm,algorithmic}
\usepackage{textcomp}
\def\BibTeX{{\rm B\kern-.05em{\sc i\kern-.025em b}\kern-.08em
    T\kern-.1667em\lower.7ex\hbox{E}\kern-.125emX}}

\usepackage{comment}
\excludecomment{toexclude}

\usepackage{xcolor}
\definecolor{UPB_blue}{HTML}{18b0e2}
\definecolor{UPB_blue_light}{HTML}{deebf7}

\usepackage{tikz}
\usetikzlibrary{shapes.multipart}
\usetikzlibrary {arrows.meta}
\usetikzlibrary{positioning}
\usetikzlibrary{fit}
\usetikzlibrary{calc}
\usetikzlibrary{decorations.markings} 
\pgfdeclarelayer{background}
\pgfsetlayers{background,main}
\colorlet{inbox}{lightgray!20}
\colorlet{outbox}{lightgray!50}

\usepackage[pdffitwindow=false,
            plainpages=false,
            pdfpagelabels=true,
            pdfpagemode=UseOutlines,
            pdfpagelayout=SinglePage,
            bookmarks=false,
            colorlinks=true,
            hyperfootnotes=false,
            linkcolor=blue,
            urlcolor=blue!30!black,
            citecolor=green!50!black]{hyperref}

\newcommand{\scs}{\mathcal{S}}

\newcommand{\tus}{{\tau_s}}

\newcommand{\im}{i=1,\ldots,m}
\newcommand{\jn}{j=1,\ldots,n}
\newcommand{\intd}{\text{d}}

\newcommand{\vareps}{\varepsilon}

\newcommand{\norm}[1]{\left\lVert#1\right\rVert}

\newtheorem{theorem}{Theorem}

\newtheorem{assumption}{Assumption}
\newtheorem{remark}{Remark}

\newcommand{\cC}{\mathcal{C}}

\newcommand{\cW}{\mathcal{W}}

\newcommand{\cK}{\mathcal{K}}

\newcommand{\cL}{\mathcal{L}}

\newcommand{\RN}[1]{\uppercase\expandafter{\romannumeral#1}}
\newcommand{\eps}{\varepsilon}

\newcommand{\N}{{\mathbb{N}}}

\newcommand{\R}{{\mathbb{R}}}

\newcommand{\setdef}[2]{\left\{\, #1 \left|\, \vphantom{#1} #2\right.\right\}}

\newcommand{\diag}{\text{\rm diag\,}}
\newcommand{\ds}[1]{{\rm \, d} #1 \,}

\DeclareMathOperator{\loc}{loc}

\newcommand{\blue}[1]{#1}
\newcommand{\bluen}[1]{#1}

\newlength{\innersep}
\newlength{\maxlength}
\newlength{\dummylength}

\newcommand{\UpperNilBlock}[3]{
\setlength{\arraycolsep}{0pt}
\renewcommand{\arraystretch}{0}
\settowidth{\maxlength}{$#1$}
\settoheight{\dummylength}{$#1$}
\ifdim\dummylength>\maxlength
  \setlength{\maxlength}{\dummylength}
\fi
\settowidth{\dummylength}{$#2$}
\ifdim\dummylength>\maxlength
  \setlength{\maxlength}{\dummylength}
\fi
\settoheight{\dummylength}{$#2$}
\ifdim\dummylength>\maxlength
  \setlength{\maxlength}{\dummylength}
\fi
\setlength{\innersep}{0.1\maxlength}
\addtolength{\maxlength}{\innersep}
\addtolength{\maxlength}{\innersep}
\newcommand{\invisiblebox}{\phantom{\rule{\maxlength}{\maxlength}}}
\begin{array}{cccc}
  \tikz[remember picture]{\node[outer sep=0,inner sep=\innersep] (a11) {$#1$};} & \tikz[remember picture]{\node[outer sep=0,inner sep=\innersep] (a21) {$#2$};}  &  & \invisiblebox\\
  &  & \phantom{\rule{#3}{#3}} &  \\
  \invisiblebox & \invisiblebox & & \tikz[remember picture]{\node[outer sep=0,inner sep=\innersep] (a43) {$#2$};}\\
   \invisiblebox & \invisiblebox &   & \tikz[remember picture]{\node[outer sep=0,inner sep=\innersep] (a44) {$#1$};}
\end{array}
\tikz[remember picture, overlay,shorten >=4pt,shorten <=4pt] \draw (a11.center) edge[very thick] (a44.center);
\tikz[remember picture, overlay,shorten >=4pt,shorten <=4pt] \draw (a21.center) edge[very thick] (a43.center);
}

\newcommand{\SmallNilBlock}[3]{
\UpperNilBlock{\mbox{\scriptsize $#1$}}{\mbox{\scriptsize $#2$}}{#3}
}

\newenvironment{smallpmatrix}
{\left(\begin{smallmatrix}}
{\end{smallmatrix}\right)}

\newenvironment{smallbmatrix}
{\left[\begin{smallmatrix}}
{\end{smallmatrix}\right]}

\DeclareRobustCommand{\svdots}{%
  \vcenter{%
    \offinterlineskip
    \hbox{.}
    \vskip0.2\normalbaselineskip
    \hbox{.}
    \vskip0.2\normalbaselineskip
    \hbox{.}%
  }%
}

\begin{document}
\title{Prescribed Performance Control of Uncertain \blue{Higher Relative Degree} Nonlinear Systems in the Presence of Delays}

\author{Thomas Berger\thanks{T. Berger acknowledges funding by the Deutsche Forschungsgemeinschaft (DFG, German
Research Foundation) -- Project-ID 471539468.}, Lampros N. Bikas, Jan Hachmeister, and George A. Rovithakis
\thanks{Thomas Berger is with the Martin Luther University Halle-Wittenberg, Institute for Mathematics, Theodor-Lieser-Str. 5, 06120 Halle (Saale), Germany (e-mail: thomas.berger@mathematik.uni-halle.de)}
\thanks{Jan Hachmeister is with the University of Paderborn, Institute for Mathematics, Warburger Str.~100, 33098~Paderborn, Germany (e-mail: janha@mail.uni-paderborn.de).}
\thanks{Lampros N. Bikas and George A. Rovithakis are with  the Department of Electrical and Computer Engineering, Aristotle University of Thessaloniki, 54124 Thessaloniki, Greece (e-mail: lnmpikas@ece.auth.gr, rovithak@ece.auth.gr).}\vspace*{-10mm}}

\maketitle

\begin{abstract}
We propose a novel feedback controller for a class of uncertain higher \blue{relative degree} nonlinear systems, subject to delays in both state measurement and control input signals. Building on the prescribed performance control framework, a delay-dependent performance correction mechanism is introduced to ensure the boundedness of all signals in the closed-loop and to keep the output tracking error strictly within a dynamically adjusted performance envelope. This mechanism adapts in response to large delays that may cause performance degradation. In the absence of delays, the correction term vanishes, and the controller recovers the nominal (user-defined) performance envelope. The effectiveness of the proposed approach is validated through simulation studies.
\end{abstract}


\section{Introduction}\label{sec:introduction}

In recent years, significant research efforts have been devoted to the design of controllers for nonlinear systems. Early works often assumed ideal communication between the controller and system, including measurement and control input signals. However, from both theoretical and practical perspectives, delays are inherent in most control systems and can significantly affect performance and stability ~\cite{Heemels2010}. Therefore, designing robust controllers that account for communication delays is of paramount importance. It is especially crucial to consider both measurement and control input delays, as these are common in practical applications, and recent progress has been reported in this direction ~\cite{Karafyllis2012,Selivanov2016,Zhou2017,Le2018,Weston2019,Battilotti2020,Nozari2020,Zhao2021,Sun2021,Yu2025}. 

A common focus of these studies has been on establishing stability conditions, often without explicit consideration of performance objectives. In trajectory tracking problems for uncertain nonlinear systems, beyond stability, it is essential to ensure strict bounds on the output tracking error, encompassing transient and steady-state behaviors. Ideally, these bounds should be prescribed and user-defined. 

Two prominent methodologies addressing this are funnel control (FC) and prescribed performance control (PPC). FC was introduced in~\cite{IlchRyan02b}, with a recent survey provided in~\cite{BergIlch21}. PPC was initially presented in~\cite{Bech2008} with further developments in~\cite{RoviSurvey}. When assuming delay-free communication, both methods can guarantee prescribed performance characteristics—such as maximum overshoot, convergence rate, and steady-state error. \blue{However, even under very small measurement delays both methods fail (cf.\ Section~\ref{sec:sim}), thus necessitating appropriate modifications.}

In the presence of measurement and control input delays, approaches such as the bang-bang controller based on FC were proposed in~\cite{Liberzon2013} where constant delays were considered. However, these methods involve increased complexity and require high-order derivatives of reference signals. Under the PPC framework, an initial attempt to address communication delays was made in~\cite{BikaRovi23}, where the control design aimed to satisfy prescribed performance objectives. However, the boundedness of all closed-loop signals is only assured for relative degree one systems under the proposed control scheme. 

In order to extend prescribed performance guarantees under delays to uncertain nonlinear systems of arbitrary relative degree, we propose a modified control scheme based on~\cite{BikaRovi23}, enhanced with a delay-dependent performance correction term. This controller enforces prescribed performance for a class of uncertain higher-order nonlinear systems under constant measurement delays and time-varying control input delays. Our approach guarantees the boundedness of all signals in the closed-loop and ensures that the output tracking error evolves strictly within a dynamically adjusted performance envelope, thereby overcoming the limitations of~\cite{BikaRovi23}. The correction term explicitly depends on the delay magnitude and vanishes in the delay-free case, thereby recovering the nominal prescribed performance. The proposed controller is of low-complexity, maintaining simplicity in implementation.


\paragraph{Notation}

In the following let $\N$ denote the natural numbers, and $\R_{\ge \tau} =[\tau,\infty)$ for $\tau\in\R$. By $\|x\|$ we denote the Euclidean norm of $x\in\R^n$. For some interval $I\subseteq\R$, some $V\subseteq\R^m$ and $k\in\N$, $\cL^\infty(I, \R^{n})$ $\big(\cL^\infty_{\rm loc} (I, \R^{n})\big)$ is the Lebesgue space of measurable, (locally) essentially bounded {functions} $f\colon I\to\R^n$, $\cW^{k,\infty}(I,  \R^{n})$ is the Sobolev space of all functions
$f:I\to\R^n$ with $k$-th order weak derivative $f^{(k)}$ and $f,f^{(1)},\ldots,f^{(k)}\in L^\infty(I, \R^{n})$, and
 $\cC^k(V,  \R^{n})$ is the set of  $k$-times continuously differentiable functions  $f:  V  \to \R^{n}$, with $\cC(V,  \R^{n}) := \cC^0(V,  \R^{n})$. Furthermore, $\cK_\infty$ denotes the set of continuous, strictly increasing and unbounded functions, and  $\cK\cL$ is the set of continuous functions, strictly decreasing with limit zero in the first argument and strictly increasing in the second argument. 

\section{Problem statement}\label{Sec:ProbState}

\paragraph{System class}

We consider nonlinear multi-input, multi-output systems  \blue{with relative degree} $n$ of the form
\begin{equation}\label{eq:Sys}
\begin{aligned}
  \dot x_{i,j}(t) &= f_{i,j}(t,\bar x_j(t)) + x_{i,j+1}(t),\\
  \dot x_{i,n}(t) &= f_{i,n}(t,\bar x_n(t),\eta(t)) \\
  &\quad + \sum_{k=1}^m g_{i,k}(t,\bar x_n(t),\eta(t)) u_k(t-\tau_u(t)),\\
  \dot \eta(t) &= h(t,\bar x_n(t),\eta(t)),
\end{aligned}
\end{equation}
for $j=1,\ldots,n-1$, $i=1,\ldots,m$, where $\bar x_j = (x_{1,1}, \ldots, x_{1,j}, \ldots, x_{m,1},\ldots, x_{m,j})^\top$ is the part of the state available for measurement, $\eta:\R_{\ge 0}\to\R^q$ is the state of the internal dynamics, $u_i$ are the control inputs and $y_i:=x_{i,1}$ are the system outputs. The functions $h,f_{i,j},g_{i,k}$, $i,k=1,\ldots,m$, are assumed to be piecewise continuous and bounded in~$t$ and locally Lipschitz in~$(\bar x_n, \eta)$ for $j=n$, or locally Lipschitz in~$\bar x_j$ for $j=1,\hdots,n-1$, respectively.

We consider the problem of output tracking control with the objective of achieving a prescribed performance of the tracking error, under the effect of time-varying input delays (described by~$\tau_u(t)\ge 0$) and constant state measurement delays (described by~$\tau_s\ge0$). 
The latter means that only the delayed information $x_{i,j}(t\!-\!\tau_s)$ of the state measurement is available for controller design. Because of this delay, in order for the problem to be well-posed, an initial history of the state is required on the interval $[-\tau_s-\bar \tau_u,0]$, where $\bar \tau_u \ge \tau_u(t)$ for all $t\ge 0$, i.e.,
\begin{equation}\label{eq:initial_history}
    (\bar x_n,\eta)|_{[-\tau_s-\bar\tau_u,0]} = \varphi = (\bar x_n^\varphi,\eta^\varphi)\in \cC([-\tau_s-\bar\tau_u,0],\R^{nm+q}).
\end{equation}
For $u_i\in \cL^\infty_{\loc}([-\tau_s-\bar\tau_u,\infty),\R)$, $\im$, we call $ (\bar x_n,\eta):[-\tau_s-\bar\tau_u,\omega)\to\R^{nm+q}$, $\omega\in(0,\infty]$, a solution of~\eqref{eq:Sys},~\eqref{eq:initial_history}, if it is locally absolutely continuous and satisfies~\eqref{eq:Sys} for almost all $t\in[0,\omega)$. A solution is called maximal, if it has no right extension that is also a solution; it is global, if $\omega=\infty$.

We need the following additional assumption on the nonlinearities $f_{i,n}$ and $h$.

\begin{assumption}\label{Ass1}
For each $i=1,\ldots,m$ there exists $d_i\in \cL^\infty(\R_{\ge 0},\R)$ such that, for all $(t,x,\eta)\in\R_{\ge 0}\times\R^{mn+q}$,
\[
     |f_{i,n}(t,x,\eta)| \le |d_i(t)| (\|\bar x_n\| + \|\eta\| + 1).
\]
Furthermore, the last equation in~\eqref{eq:Sys} is practically input-to-state stable in the sense that there exists a constant $c\ge 0$, a $\cK\cL$-function $\kappa$ and a $\cK_\infty$-function $\gamma$ such that, for any $\xi\in \cL^\infty_{\loc}(\R_{\ge 0}, \R^{nm})$ and any $\eta^0\in\R^q$ the solution of $\dot \eta(t) = h(t,\xi(t),\eta(t))$, $\eta(0)=\eta^0$, is global and satisfies
\[
    \forall\, t\ge 0:\ \|\eta(t)\| \le \kappa(t,\|\eta^0\|) + \gamma\big(\sup\nolimits_{s\in[0,t]} \|\xi(s)\|\big) + c.
\]
Additionally, the function~$\gamma$ is linearly bounded, i.e., $\gamma(x)\le \bar \gamma_1 x + \bar\gamma_2$ for all $x\ge 0$ and some $\bar\gamma_1,\bar\gamma_2 >0$. 
\end{assumption}


\begin{remark}
    We note that the Lipschitz assumption on the nonlinearities~$f_{i,n}$ cannot be waived in general. For instance, consider the system $\dot x(t) = x(t)^2 + u(t-\tau)$, $x|_{[-\tau,0]} \equiv x^0> 0$.
    Since the only reasonable choice for the input for $t\in [-\tau,0]$ is $u(t)=0$, as the system is not ``active'' yet, this leads to the initial-value problem $\dot x(t) = x(t)^2$, $x(0)=x^0$, $t\in[0,\tau]$,
    the solution of which is given by $x(t) = \left(\frac{1}{x^0} - t\right)^{-1}$, $t\in[0,\min\{\tau,1/x^0\})$.
    If $\frac{1}{x^0} < \tau$, then this leads to a blow-up of the solutions, hence existence of global solutions of the closed-loop system cannot be guaranteed by any control algorithm.
\end{remark}

\begin{assumption}\label{Ass2}
The delays $\tau_u\in \cC^1(\R_{\ge 0},\R)$ and $\tau_s\ge 0$ are assumed to be known exactly and bounded such that $\bar \tau_u \ge \tau_u(t)$ for all $t\ge 0$. Furthermore, we assume that there exists $\dot{\bar{\tau}}_u<1$ such that $\dot \tau_u(t) \le \dot{\bar{\tau}}_u$ for all $t\ge 0$. 
\end{assumption} 

The strict boundedness requirement on $\dot{\tau}_u(t)$ is necessary to guarantee the satisfaction of the first-in first-out principle~\cite{BikaRovi23}. The required knowledge of the delays $\tau_s(\cdot)$ and $\tau_u(\cdot)$ might seem to be a strong assumption, but in many applications they can indeed be estimated very well, or, by intentionally delaying some measurements, a prescribed quantity can be achieved. 


\begin{toexclude}

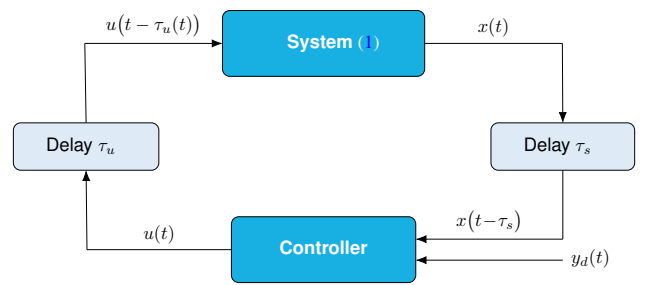
\begin{figure}[h!]
\resizebox{0.5\textwidth}{!}{
\begin{tikzpicture}[every text node part/.style={align=center},
rectanglenode1/.style={rectangle, rounded corners, draw=black, fill=UPB_blue, ultra thin, text=white},
rectanglenode2/.style={rectangle, rounded corners, draw=black, fill=UPB_blue_light, ultra thin, minimum size=5mm},
]

\textsf{
\node[rectangle]	      				(middle)       		 						{$\qquad$\hspace{100pt} \\ $\qquad$};
\node[rectanglenode1]      				(system)      		[above=of middle]		{\quad\\[5pt] $\qquad$ \textbf{\, \, System~\eqref{eq:Sys}} $\qquad$\! \\ \quad};
\node[rectanglenode1]      				(controller)		[below=of middle] 		{\quad\\[5pt] $\qquad$ \textbf{Controller} $\qquad$ \\ \quad};
\node[rectanglenode2]        			(delay_s)       	[right=of middle,xshift=-0.5mm] 		{\quad\\[0pt] $\quad\, $ Delay $\tau_s \qquad$ \\[-5pt] \quad};
\node[rectanglenode2]        			(delay_u)			[left=of middle,xshift=1mm] 		{\quad\\[0pt] $\quad\, $ Delay $\tau_u  \qquad$ \\[-5pt] \quad};
}
\node[rectangle] (upper_right) [right=of system] 		{\quad $\hspace{75pt}$};
\node[rectangle] (upper_left)  [left=of system] 		{\quad $\hspace{75pt}$};
\node[rectangle] (lower_right) [right=of controller,xshift=-0.5mm] 	{\quad $\hspace{88pt}$ \\[5pt] $\hspace{30pt} y_d(t)$};
\node[rectangle] (lower_left)  [left=of controller,xshift=1mm]	{\quad $\hspace{88pt}$};

\draw (system.east) -- (upper_right.center) node[above,midway] {$x(t)$};
\draw[-{Latex[length=2mm]}] (upper_right.center) -- (delay_s.north);

\draw (delay_u.north) -- (upper_left.center);
\draw[-{Latex[length=2mm]}] (upper_left.center) -- (system.west) node[above,midway] {$u\big(t-\tau_u(t)\big)$};

\draw (delay_s.south) -- ([yshift=2mm]lower_right.center);
\draw[-{Latex[length=2mm]}] ([yshift=2mm]lower_right.center) -- ([yshift=2mm]controller.east) node[above,midway] {$x\big(t\!-\!\tau_s\big)$};

\draw (controller.west) -- (lower_left.center) node[above,midway]{$u(t)$};
\draw[-{Latex[length=2mm]}] (lower_left.center) -- (delay_u.south);

\draw[-{Latex[length=2mm]}] ([yshift=-2mm]lower_right.center) -- ([yshift=-2mm]controller.east);
\end{tikzpicture}
}
\caption{Structure of the closed-loop system.}
\label{Fig:1}
\end{figure}
\end{toexclude}

\blue{
\begin{remark}\label{Rem:dist_u}
    Note that it is possible to incorporate the occurrence of input disturbances in the system class, i.e., $u_k \to u_k + w_k$ for some $w_k\in \cL^\infty([-\tau_s-\bar\tau_u,\infty),\R)$ in~\eqref{eq:Sys}. In this case, the functions $f_{i,n}$ in the second of equations~\eqref{eq:Sys} can be rewritten as
    \[
        f_{i,n}(t,\bar x_n,\eta) \to f_{i,n}(t,\bar x_n,\eta) + \sum_{k=1}^m g_{i,k}(t,\bar x_n,\eta) w_k(t-\tau_u(t)).
    \]
    If all $g_{i,k}$ satisfy Assumption~\ref{Ass1}, then also the above remodeled nonlinearities satisfy  Assumption~\ref{Ass1}.
\end{remark}
}

\paragraph{Control objective}

The objective is, for given functions $\psi_{i,1}\in \cW^{1,\infty}(\R_{\ge -\tau_s},\R)$ which are bounded and satisfy $\inf_{t\ge -\tau_s} \psi_{i,1}(t) >0$, and reference signals $y_{d,i}\in \mathcal{W}^{1,\infty}(\R_{\ge -\tau_s},\R)$, $i=1,\ldots,m$, to design a controller which achieves that $|y_i(t)-y_{d,i}(t)| < \psi_{i,1}(t)$ for all $i=1,\ldots,m$ and all $t\ge 0$. 
Furthermore, all closed-loop signals should remain bounded. The controller should not require any knowledge of the system parameters (initial history~$\varphi$, nonlinearities~$f_{i,j}, g_{i,k}, h$) and should be of low complexity (no approximation or adaptive structures are used to obtain that knowledge, no hard calculations are performed to create the control signal). A controller which satisfies these requirements is inherently robust with respect to uncertainties or disturbances (with the exception of measurement noise).


\begin{remark}
    The performance functions $\psi_{i,1}$, $\im$, are user-defined and can be constructed to introduce performance bounds on the output tracking error with respect to transient and steady-state behavior. A candidate selection is the exponentially decreasing $\psi_{i,1}(t)=(\lambda_{i,1}^0 - \lambda_{i,1}^\infty) e^{-c_{i,1} t} + \lambda_{i,1}^\infty$ with $\lambda_{i,1}^0>|y_i(0)-y_{d,i}(0)|$ and \blue{suitable $\lambda_{i,1}^\infty, c_{i,1}>0$.} 
\end{remark}

\begin{remark}
    The control objective, as formulated above, cannot be achieved by any control algorithm, when the control input is ``switched on'' at a certain time (usually $t=0$), and zero before (that is, $u(t) = 0$ for $t\le 0$). Nevertheless, the new design that we propose in Section~\ref{Sec:ContrDes} achieves the prescribed performance objective with an additional correction term, the magnitude of which depends on the magnitude of the delays (and vanishes for zero delay). 
\end{remark}

\section{Controller design}\label{Sec:ContrDes}
In this Section we propose a control scheme to achieve the control objectives stated in Section \ref{Sec:ProbState}. The control design philosophy is based on the line of analysis of PPC structure and incorporates a delay-dependent transformation of the output tracking error, resulting in a closed-loop system that includes at least one delay-free control input.

Let $\psi_{i,1}\in \cW^{1,\infty}(\R_{\ge -\tau_s},\R)$, $\im$, be given and select the auxiliary functions $\psi_{i,j}\in \cW^{1,\infty}(\R_{\ge -\tau_s},\R)$ such that $\inf_{t\ge 0} \psi_{i,j}(t) >-\tau_s$, for $i=1,\ldots,m$ and $j=1,\ldots,n$. Further, choose the controller design parameters:

\begin{itemize}
  \item a matrix $S = (s_{i,j}) \in\R^{m\times m}$, which is sign definite, i.e., there exist $s^*>0$ and $\sigma\in\{-1,+1\}$ such that $\sigma v^\top S v \ge s^* \|v\|^2$ for all $v\in\R^m$,
  \item some freely selected control gains $\alpha > 0$ and $k_{i,j}, k_n>0$, $i=1,\ldots,m$, $j=1,\ldots,n-1$.
\end{itemize}
The following expressions define the controller we propose for $i=1,\ldots,m$, and $t\ge 0$:
\begin{equation}\label{eq:controller}
\small{
\begin{aligned}
z_{i,1}(t) &= \frac{x_{i,1}(t\!-\!\tau_s) - y_{d,i}(t\!-\!\tau_s) + I_{i,1}(t)}{\psi_{i,1}(t\!-\!\tau_s)},\\
z_{i,j}(t) &= \frac{x_{i,j}(t\!-\!\tau_s) - a_{i,j-1}(t) \!+\! \sum\limits_{k=1}^j \tbinom{j-1}{j-k} (-\alpha)^{j-k} I_{i,k}(t)}{\psi_{i,j}(t\!-\!\tau_s)},\\
&\quad\ j=2,\ldots,n,\\
a_{i,j}(t) &= -k_{i,j} \frac{z_{i,j}(t)}{1 - (z_{i,j}(t))^2},\quad j=1,\ldots,n-1\\
z_n(t) &= (z_{1,n}(t),\ldots,z_{m,n}(t))^\top,\\
u_i(t) &= -\sigma k_n \chi(\|z_{n}(t)\|)  \frac{z_{i,n}(t)}{1 - \|z_{n}(t)\|^2},
\end{aligned}}
\end{equation}
with
\begin{equation}\label{eq:controller-I}
\small{
\begin{aligned}
\dot I_{i,j}(t)&= I_{i,j+1}(t) - \alpha I_{i,j}(t),\quad I_{i,j}(0)=0,\\
&\qquad\qquad\qquad\qquad\qquad\quad j=1,\ldots,n-1,\\
\dot I_{i,n}(t) &=  \sum_{k=1}^m s_{i,k} \big(u_k(t) - u_k(t\!-\!\tau_s\!-\!\tau_u(t\!-\!\tus))\big)\\
&\quad -\alpha I_{i,n}(t),\quad  I_{i,n}(0)=0,
\end{aligned}
}
\end{equation}
and  an activation function $\chi:[0,1] \to [0,1\!-\!\delta]$ \blue{with}
\begin{equation}\label{eq:act_fct}
    \blue{\chi(s) = 0\ \text{for}\ s\le \delta,\ \chi(s)=
        s\!-\!\delta\ \text{for}\ s >\delta,\ \ \delta\in [0,1).}
\end{equation}
Furthermore, we assume that $u_i(t) = 0$ for $t\in [-\tau_s-\bar\tau_u,0]$ and, for simplicity,
\begin{equation}\label{eq:psi-n}
    \forall\, t\ge -\tau_s:\ \psi_{1,n}(t) = \ldots = \psi_{m,n}(t) =: \psi_n(t).
\end{equation}

\blue{We like to note that the terms $I_{i,1}$ serve as \textit{correction terms} of the reference signals $y_{d,i}$, with the design philosophy being as follows: In the presence of delays, the traditional PPC law is no longer able to react to the instantaneous reference signal. To address this issue, the novel control scheme~\eqref{eq:controller} achieves that the corrected error $x_{i,1}(t\!-\!\tau_s) - y_{d,i}(t\!-\!\tau_s) + I_{i,1}(t)$ stays inside the prescribed error margin. Because of the higher order of the system, successive correction terms are necessary for the higher order terms, so that $|z_{i,j}(t)|<1$ is achieved for $i=1,\ldots,m$, $j=1,\ldots,n$. If no delays are present, then all correction terms vanish by construction, that is $I_{i,j}(t)= 0$. In the presence of delays, the correction terms are bounded as shown in Theorem~\ref{thm:FC}.} Compared to~\cite{BikaRovi23}, the correction terms are the crucial modification of the controller design, which facilitates its feasibility.

\begin{remark}
By incorporating the activation function $\chi$ to the control input, we enforce $u_i$ to be zero when $z_{n}$ is small (i.e., $\|z_{n}\|<\delta$). In such case the $I_{i,j}$-terms will converge to zero, thus recovering the original shape of the performance envelope. This modification allows the performance envelope to be adjusted only when the error evolves sufficiently close to the boundaries; a beneficial property in practice as one cannot exclude the case where despite the presence of a large delay the error evolves close to zero. In the context of funnel control and funnel MPC, activation functions have been employed in~\cite{BergDenn24,LanzDenn24}.
\end{remark}

\section{Main results}\label{sec:mainresults}

Before we state the main result of this article we require an additional assumption for feasibility of the control. Essentially, the assumption states that the delays $\tau_s$ and $\tau_u(t)$ need to be sufficiently small and, at the same time, $S$ needs to be a sufficiently good estimate of the control input matrix $G(t,x,\eta) = (g_{i,j}(t,x,\eta)) \in\R^{m\times m}$ for $(t,x,\eta)\in\R_{\ge 0}\times \R^{nm+q}$.

\begin{assumption}\label{Ass3}
Assume that
\[
    \sup_{(t,x,\eta)\in\R_{\ge 0}\times\R^{nm+q}} \|G(t,x,\eta)-S\| + \blue{(1-\delta)} C<  \blue{\tfrac{1-\delta}{2}} s^*,
\]
where 
\bluen{$C := \tilde c\,\|S\|\,C_\tau$ with
\[
  C_\tau := \tfrac{M}{\mu}\left(\tfrac{\bar{\dot\tau}_u}{1-\bar{\dot\tau}_u}
  + \Big(\mu+\tfrac{\|A\|M}{1-\bar{\dot\tau}_u}\Big)(\tau_s+\bar\tau_u)\right),
\]
}
$A = \diag(\tilde A,\ldots,\tilde A)\in\mathbb{R}^{nm\times nm}$,
$
    \tilde A = \left[\SmallNilBlock{\text{-}\alpha}{1}{2ex}\right] \in\R^{n\times n},
$
$M, \mu >0$ such that (note that $\mu < \alpha$ for $n>1$)
\begin{equation}\label{eq:est-eAt}
\forall\, t\ge 0:\ \|e^{At}\| \le M e^{-\mu t},
\end{equation}
and $\tilde c := \bluen{(1+\bar \gamma_1)} \max\left\{\alpha^{n-1},1\right\}\tbinom{n}{\lfloor\frac{n}{2}\rfloor} n\cdot \max_{\im}\|d_i\|_\infty +\max\left\{\alpha^{n},1\right\}\tbinom{n}{\lfloor\frac{n}{2}\rfloor}$.
\end{assumption}

\begin{remark}
    We like to note that Assumption~\ref{Ass3} is indeed quite restrictive and requires a significant amount of knowledge of the control input matrix~$G(\cdot)$, only up to slight uncertainties. Furthermore, since~$S$ is a constant matrix, $G(\cdot)$ is restricted to be globally bounded. In future research, we aim to relax these assumptions.
\end{remark}

The following result shows feasibility of the application of the controller~\eqref{eq:controller},~\eqref{eq:controller-I} to system~\eqref{eq:Sys} in the presence of delays.

\begin{theorem}\label{thm:FC}
  Consider the system~\eqref{eq:Sys} satisfying Assumptions 1 and 2. Further consider the controller~\eqref{eq:controller},~\eqref{eq:controller-I} with reference signals $y_{d,i}\in \mathcal{W}^{1,\infty}(\R_{\ge -\tau_s},\R)$, $\im$, and design parameters satisfying Assumption 3. Let $\varphi$ as in~\eqref{eq:initial_history} be an initial history such that all signals in~\eqref{eq:controller},~\eqref{eq:controller-I} are well-defined for $t\in [-\bar\tau_u,\tau_s]$ (with $I_{i,j}(t)=0$ for $t\le 0$) and, in particular, $|z_{i,j}(t)| < 1$ for $i=1,\ldots,m$, $j=1,\ldots,n-1$ and $\|z_n(t)\| < 1$. Then the application of the controller~\eqref{eq:controller},~\eqref{eq:controller-I} to the system~\eqref{eq:Sys} leads to a closed-loop system,  which has a solution and every solution can be  extended to a maximal solution $(\bar x_n, \eta, I_{1,1},\ldots, I_{m,n}):[-\tau_s-\bar \tau_u,\omega) \to \R^{2mn+q}$, $\omega\in (0,\infty]$, with the properties:
  \begin{enumerate}
      \item[(i)] global existence: $\omega = \infty$;
      \item[(ii)] all closed-loop signals $x_{i,j}, I_{i,j}, \eta$, $i=1,\ldots,m$, $j=1,\ldots,n$, and the control input $u$ are bounded;
      \item[(iii)] the system output exhibits a prescribed performance in the sense that, for all $t\ge 0$ and  $i=1,\ldots,m$,
      \[
    |y_{i}(t\!-\!\tau_s) - y_{d,i}(t\!-\!\tau_s) + I_{i,1}(t)| < \psi_{i,1}(t\!-\!\tau_s);
  \]
\item[(iv)] \bluen{the deviation from the nominal performance envelope is quantified
by the correction terms: for all $t\ge 0$ and $i=1,\ldots,m$,
\[
  |I_{i,1}(t)| \le \|S\|\,C_\tau\sup_{s\in[0,t]}\|u(s)\|,
  \ \ u:=(u_1,\ldots,u_m)^\top,
\]
with $C_\tau$ as in Assumption~3.} 
  \end{enumerate}
\end{theorem}

The proof of Theorem~\ref{thm:FC} is relegated to the Appendix.

\section{Simulations}\label{sec:sim}
\begin{figure}[bt]
\begin{center}
\includegraphics[width=0.3\textwidth,keepaspectratio]{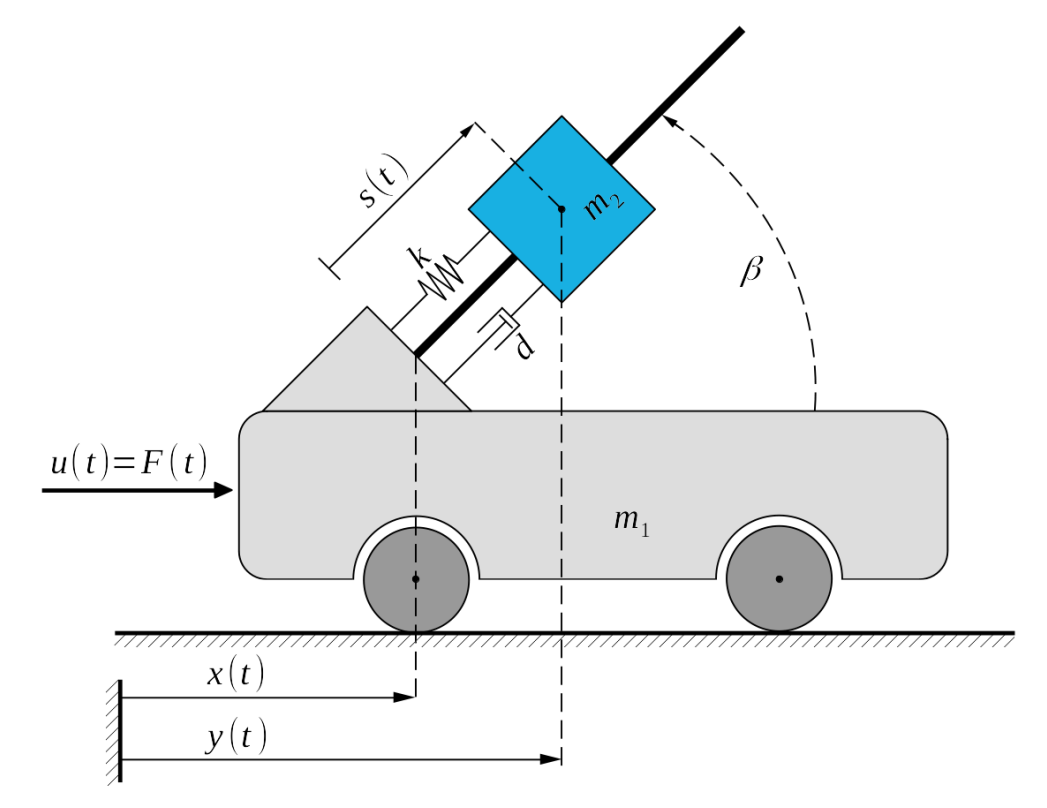}
\caption{Mass on car system}\vspace*{-8mm}
\label{Fig:2}
\end{center}
\end{figure}

To illustrate the applicability of our controller we consider a mass-spring system mounted on a moving car, see Fig. \ref{Fig:2}. This system was originally presented in \cite{SeifBlaj13} and has been used to demonstrate the application of a funnel controller in a delay-free setting in \cite{BergLeReis2018}.

The car with mass $m_1$ (in $kg$) moves horizontally and is actuated by input force $u(t)=F(t)$ (in $N$). On the car, the mass $m_2$ (in $kg$) is mounted via a spring-damper combination and moves along an axis inclined by the angle $\beta$ (in $rad$). The equations of motion are
\begin{equation}\label{4.36}
\begin{smallbmatrix}m_1+m_2 & m_2\cos\beta \\ m_2\cos\beta & m_2\end{smallbmatrix}
\begin{smallpmatrix}\ddot x(t) \\ \ddot s(t)\end{smallpmatrix} +
\begin{smallpmatrix}0 \\ k s(t)+d \dot s(t)\end{smallpmatrix} =
\begin{smallpmatrix}u(t) \\ 0\end{smallpmatrix}
\end{equation}
where $x(t)$ is the horizontal position of the car and $s(t)$ is the relative position of mass $m_1$ on the car along the inclined axis. The parameters $c$ (in $N/m$) and $d$ (in $Ns/m$) denote the spring and damper coefficients, respectively. The system's output is the horizontal position of mass $m_1$, $y(t)=x(t)+s(t)\cos\beta$.

As we can only control the force acting on the car and not $m_1$ itself, the system possesses internal dynamics. For the simulations we choose the same parameters as in \cite{BergLeReis2018}, i.e., $m_1=4$, $m_2=1$, $k=2$ and $d=1$. The initial history of the system for $t\in[-\tau_s-\bar\tau_u,0)$ is given by $x(t)=s(t)\equiv0$ and $\dot x(t)=\dot s(t)\equiv 0$. The goal is to design a control input $u(t)$ so that the output $y(t)$ tracks the reference trajectory $y_{d}(t)=\cos t$. \blue{Both the control input and the state measurement are subject to constant delays $\tau_s$ and $\tau_u$. We conducted simulations for the case of $\beta = \frac \pi 4$, i.e., relative degree two. We demonstrate the operation of the controller under small and large delays to illustrate the delay-impact on the performance.}

\begin{figure}[bt]
\begin{center}
\includegraphics[width=0.4\textwidth,keepaspectratio]{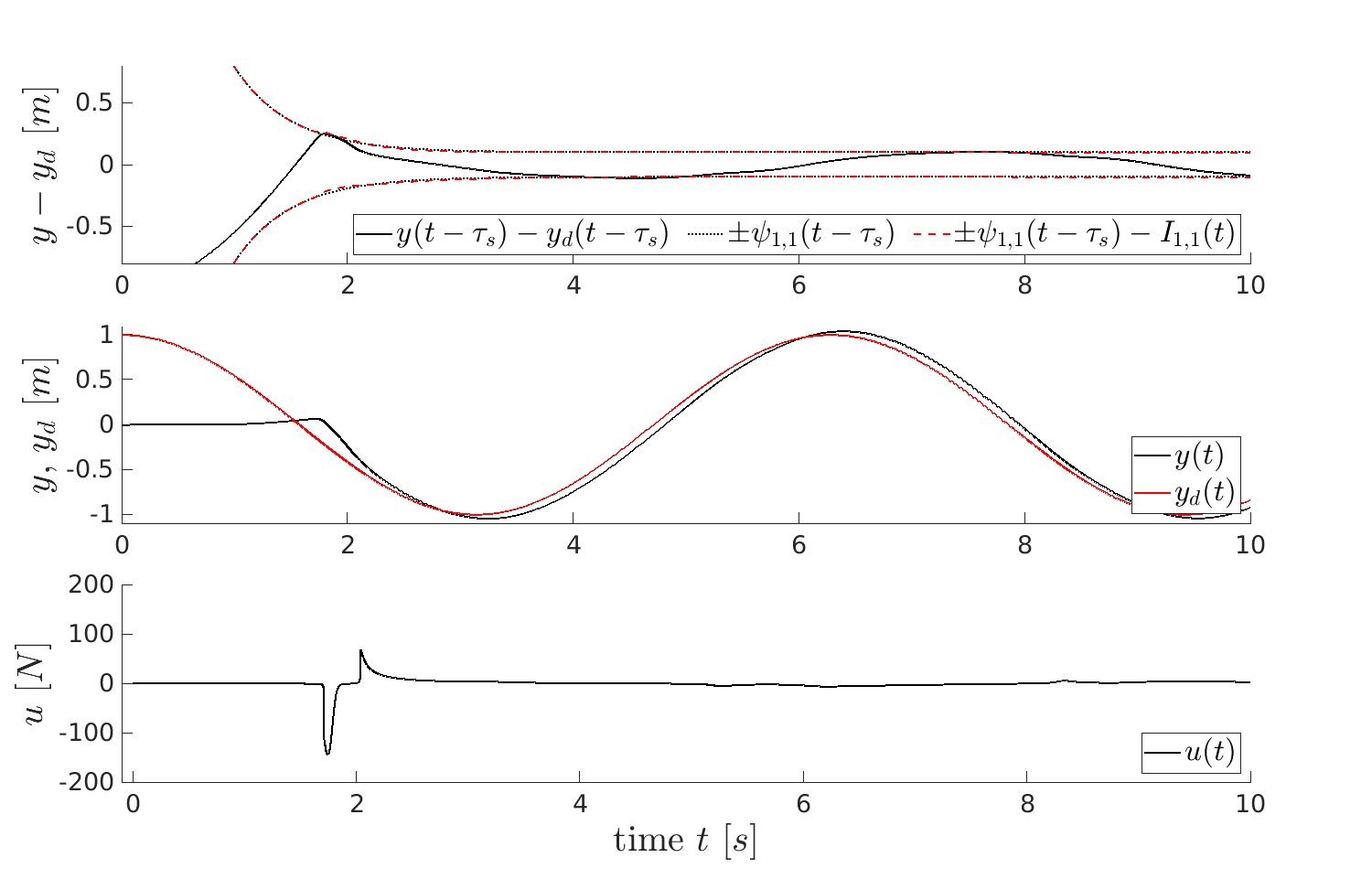}
\caption{Simulation of system~\eqref{4.36} with delays $\tau_s=\tau_u=0.01 \ \mathrm{[s]}$ under controller~\eqref{eq:controller},~\eqref{eq:controller-I} for $\beta=\frac{\pi}{4}$.}\vspace*{-5mm}
\label{Fig:3}
\end{center}
\end{figure}

\begin{figure}[bt]
\begin{center}
\includegraphics[width=0.4\textwidth,keepaspectratio]{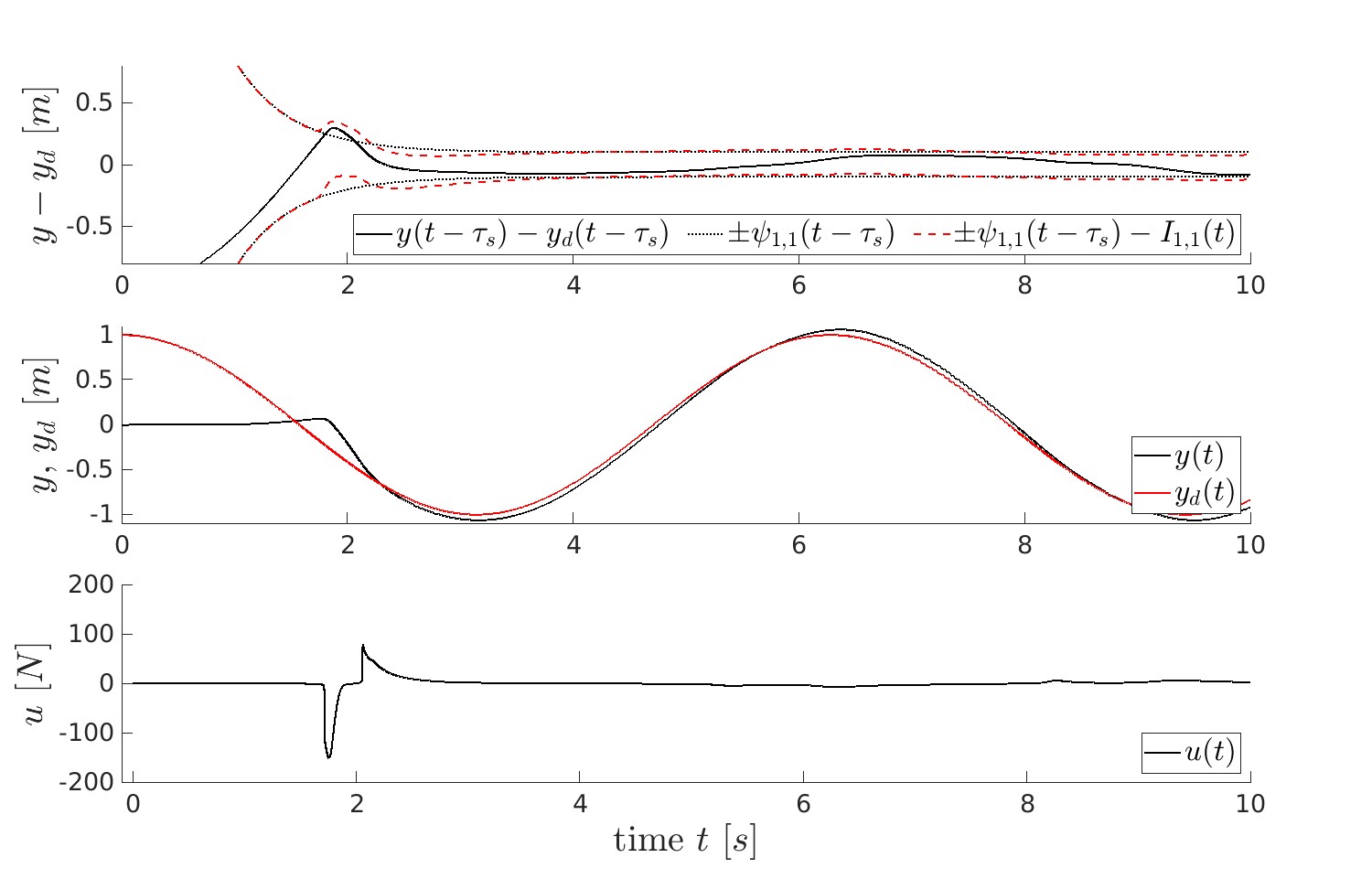}
\caption{Simulation of system~\eqref{4.36} with delays $\tau_s=\tau_u=0.05 \ \mathrm{[s]}$ under controller~\eqref{eq:controller},~\eqref{eq:controller-I} for $\beta=\frac{\pi}{4}$.}\vspace*{-5mm}
\label{Fig:31}
\end{center}
\end{figure}

\begin{figure}[bt]
\begin{center}
\includegraphics[width=0.4\textwidth,keepaspectratio]{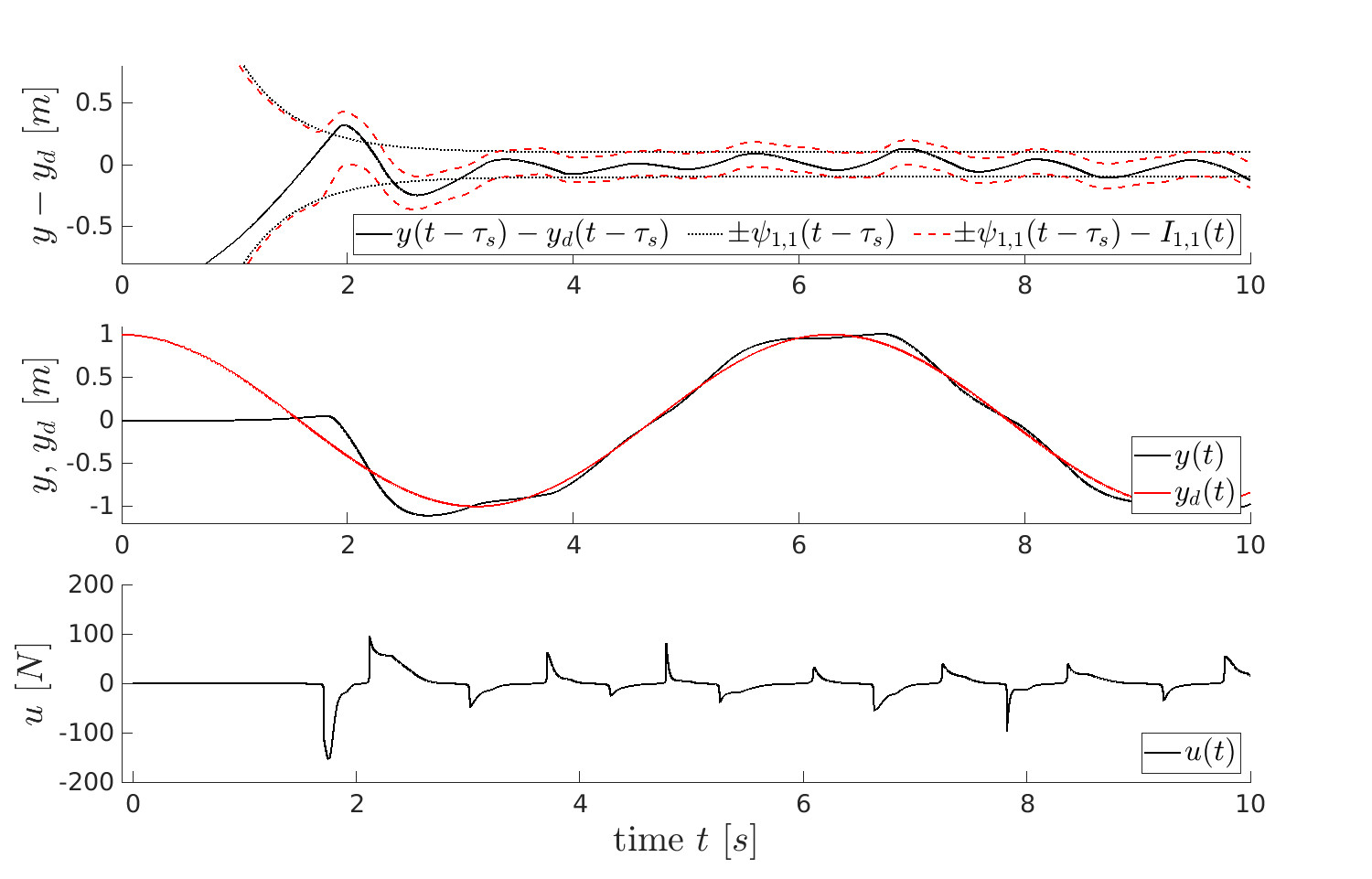}
\caption{Simulation of system~\eqref{4.36} with delays $\tau_s=\tau_u=0.1 \ \mathrm{[s]}$ and external disturbance $w(t)$ under controller~\eqref{eq:controller},~\eqref{eq:controller-I} for $\beta=\frac{\pi}{4}$.}\vspace*{-8mm}
\label{Fig:32}
\end{center}
\end{figure}

\blue{For $\beta = \frac \pi 4$} the system is of the form~\eqref{eq:Sys} with $n=q=2$ and $m=1$ as well as $G(t,x,\eta)=\tfrac19$ as shown in~\cite{BergIlch21}. We choose the controller design parameters $S=\frac{1}{9}$, $k_{1,1} = k_2 =1$ and $\alpha = 1$. The activation function is chosen as $\chi(s) = s$ with $\delta=0$. 
For the performance functions we choose $\psi_{1,1}(t)=5e^{-2t}+0.1$ and $\psi_{1,2}(t)=10e^{-2t}+0.5$. The application of the controller~\eqref{eq:controller},~\eqref{eq:controller-I} to~\eqref{4.36} with delays $\tau_s=\tau_u=0.01 \ \mathrm{[s]}$ is depicted in Fig.~\ref{Fig:3}, showing the output tracking error alongside the funnel boundary in the top image, the output alongside the reference signal in the middle and the control input signal in the bottom image. We see that apart from some limited time intervals, the error $y(t-\tau_s)-y_d(t-\tau_s)$ remained within the original funnel $\pm \psi_{1,1}(t-\tau_s)$. \blue{Subsequently, we increased the delays to $\tau_s=\tau_u=0.05 \ \mathrm{[s]}$ and $\tau_s=\tau_u=0.1 \ \mathrm{[s]}$. For the latter case we also incorporated bounded, additive external disturbances $w(t)$ in the system dynamics (i.e., $u+w$ in the right hand side of~\eqref{4.36}) with $w(t)=5\sin(5t)$, $t\geq-\tau_s-\bar \tau_u$, which is, in view of Remark~\ref{Rem:dist_u}, feasible since $G(\cdot)$ is bounded. The results are depicted in Fig.~\ref{Fig:31} and Fig.~\ref{Fig:32}, respectively. It is clearly shown that larger delays, combined with the external disturbance, led to larger $I_{1,1}$-terms, in accordance with Theorem~\ref{thm:FC}\,(iv), and consequently, to increased error bounds.}

As mentioned before, the same system has already been considered in a delay-free setting in~\cite{BergLeReis2018}. However, when using the controller therein and applying a state measurement delay~$\tau_s$ as little as $0.005$ [$s$] and no control input delay~$\tau_u$, the closed-loop system becomes unstable and the simulation fails at $t\approx 2.5$, as the controller is not able to keep the error between the funnel boundaries.

\newpage



\section{Conclusion}\label{sec:conclusion}
A novel control approach is presented for a class of uncertain nonlinear systems of arbitrary relative degree, addressing the challenges of state and control input delays. Extending the PPC framework, we introduce a delay-dependent correction term that actively compensates for communication delays. This results in guaranteed, dynamically adjusted output tracking performance. A key advantage of our controller is its adaptability: it automatically recovers nominal performance when delays are minimal, eliminating the need for re-tuning. The effectiveness of this approach is validated through simulations.


\section*{Appendix - Proof of Theorem~\ref{thm:FC}}

The proof of the theorem consists of four phases. In Phase~1 we derive the differential equations of the transformed closed-loop system and in Phase~2 we guarantee the existence of solution in a maximal time interval. Further, in Phase~3 we prove that the solution evolves strictly within the prescribed performance envelope during this interval. This enables us to show that the solution is global and all signals in the closed-loop system are bounded in Phase~4.

\emph{Phase 1}: First, let $\delta_I>0$ be a constant, chosen large enough with lower bounds to be specified later in the proof. Define the set $\Omega_I:=(-\delta_I,\delta_I)$. We derive differential equations for $z_{i,j}(t)$, $\im$, $\jn$ for $t\geq \tau_s$. 
We compute
\begin{align*}
&\dot z_{i,1}(t)
=\tfrac{1}{\psi_{i,1}(t\!-\!\tus)}\Big(f_{i,1}\big(t, \bar x_1(t\!-\!\tus)\big) +  z_{i,2}(t)\psi_{i,2}(t\!-\!\tus) \\
&\quad  -\tfrac{k_{i,1}z_{i,1}(t)}{1-z_{i,1}(t)^2} -\dot y_{d,i}(t\!-\!\tus) -\dot \psi_{i,1}(t\!-\!\tus)z_{i,1}(t)\Big) \\
&\quad =:h_{i,1}\big(t,z_{1,1}(t),\ldots,z_{m,1}(t), z_{i,2}(t), I_{i,1}(t),\ldots, I_{m,1}(t)\big),
\end{align*}
where we used that we can express $\bar x_1(t\!-\!\tus)$ (and hence the dependency of $f_{i,1}$ on it) in terms of
\begin{itemize}
\item $\psi_{1,1}(t\!-\!\tus)z_{1,1}(t), \ \ldots, \ \psi_{m,1}(t\!-\!\tus)z_{m,1}(t)$,
\item $y_{d,1}(t\!-\!\tus), \ \ldots, \ y_{d,m}(t\!-\!\tus)$ and
\item $I_{1,1}(t), \ \ldots, \ I_{m,1}(t)$,
\end{itemize}
hence $h_{i,1}: \R_{\geq \tau_s}\times(-1,1)^{m+1}\times\Omega_I^{m}\to\R$ is well-defined and continuous. For $j=2,\ldots,n-1$ we obtain in a similar way that
\begin{align*}
\dot z_{i,j}(t)
&=\tfrac{1}{\psi_{i,j}(t\!-\!\tus)}\Big(z_{i,j+1}(t)\psi_{i,j+1}(t\!-\!\tus)-\tfrac{k_{i,j}z_{i,j}(t)}{1-z_{i,j}(t)^2}\\
&\quad + f_{i,j}(t,\bar x_j(t\!-\!\tus)) -\dot \psi_{i,j}(t\!-\!\tus)z_{i,j}(t)-\dot a_{i,j-1}(t)\Big)\\
&=:h_{i,j}\big(t,z_{1,1}(t),\ldots,z_{m,1}(t),\ldots,z_{1,j}(t),\ldots, z_{m,j}(t),\\
&\quad z_{i,j+1}(t),I_{1,1}(t),\ldots,I_{m,1}(t),\ldots,I_{1,j}(t),\ldots,I_{m,j}(t)\big)
\end{align*}
where we used that
\begin{align*}
&\sum\limits_{k=1}^j
\tbinom{j-1}{j-k}\Big((-\alpha)^{j+1-k}I_{i,k}(t)+(-\alpha)^{j-k}I_{i,k+1}(t)\Big)\\
&= \sum\limits_{k=1}^{j+1}\tbinom{j}{j+1-k}(-\alpha)^{j+1-k}I_{i,k}(t).
\end{align*}
Furthermore,  we can express $\bar x_j(t\!-\!\tus)$ (and hence the dependency of $f_{i,j}$ on it) in terms of
\begin{itemize}
\item $\psi_{1,1}(t\!-\!\tus)z_{1,1}(t), \ \ldots, \ \psi_{m,j}(t\!-\!\tus)z_{m,j}(t)$,
\item $y_{d,1}(t\!-\!\tus), \ \ldots, \ y_{d,m}(t\!-\!\tus)$,
\item $a_{1,1}(t), \ \ldots, \ a_{m,j-1}(t)$ and
\item $I_{1,1}(t), \ \ldots, \ I_{m,j}(t)$,
\end{itemize}
thus $h_{i,j}: \R_{\geq \tau_s}\times(-1,1)^{mj+1}\times\Omega_I^{mj}\to\R$ is well-defined and continuous. For $j=n$ we obtain
\begin{align*}
&\dot z_{i,n}(t)
=\frac{1}{\psi_n(t\!-\!\tus)}\Bigg(f_{i,n}(t\!-\!\tus,\bar x_n(t\!-\!\tus),\eta(t\!-\!\tus))-\dot a_{i,n-1}(t)\\
&\quad +\sum\limits_{k=1}^m\big(g_{i,k}(t\!-\!\tus,\bar x_n(t\!-\!\tus),\eta(t\!-\!\tus))\\
&\qquad\qquad - s_{i,k}\big) u_k(t\!-\!\tus\!-\!\tau_u(t\!-\!\tus))\\
&\quad-\sum\limits_{k=1}^{n}\tbinom{n}{n+1-k}(-\alpha)^{n+1-k}I_{i,k}(t)\\
&\quad+\sum\limits_{k=1}^m s_{i,k} u_k(t) -\dot \psi_n(t\!-\!\tus)z_{i,n}(t)\Bigg)\\
&=:h_{i,n}\big(t,z_{1,1}(t),\ldots,z_{m,1}(t),\ldots,z_{1,n-1}(t),\ldots, z_{m,n-1}(t),\\
&\quad z_n(t),z_{n}(t\!-\!\tus\!-\!\tau_u(t\!-\!\tus)),I_{1,1}(t),\ldots,I_{m,n}(t),\eta(t\!-\!\tus)\big)
\end{align*}
where we can express $\bar x_n(t\!-\!\tus)$ (and hence the dependency of $f_{i,n}$ and $g_{i,k}$ on it) in terms of
\begin{itemize}
\item $\psi_{1,1}(t\!-\!\tus)z_{1,1}(t), \ \ldots, \ \psi_{m,n}(t\!-\!\tus)z_{m,n}(t)$,
\item $y_{d,1}(t\!-\!\tus), \ \ldots, \ y_{d,m}(t\!-\!\tus)$,
\item $a_{1,1}(t), \ \ldots, \ a_{m,n-1}(t)$ and
\item $I_{1,1}(t), \ \ldots, \ I_{m,n}(t)$,
\end{itemize}
i.e., there exists a continuous function $F_n: \R_{\geq \tau_s}\times(-1,1)^{m(n-1)}\times \Omega \times\Omega_I^{mn}\to\R^{mn}$,
such that, for all $t\ge \tau_s$, 
\begin{multline*}
    \bar x_n(t\!-\!\tus) = F_n\big(t,z_{1,1}(t),\ldots,z_{m,1}(t),\ldots,z_{1,n-1}(t),\\
    \ldots, z_{m,n-1}(t), z_n(t), I_{1,1}(t),\ldots,I_{m,n}(t)\big).
\end{multline*}
The right-hand side of the differential equation for $z_{i,n}$ is then a continuous function $h_{i,n}: \R_{\geq \tau_s}\times(-1,1)^{m(n-1)}\times\Omega^2\times\times \Omega_I^{mn} \times \R^{q}\to\R$,
where $\Omega:=B(0,1)\subset\R^m$ is the open unit ball in $\R^m$.

\emph{Phase 2}: We show existence of a local solution. Considering~\eqref{eq:controller-I} we have
\begin{equation}\label{eq:ODE_Iij}
\begin{smallpmatrix}\dot I_{1,1}(t) \\ \svdots \\ \dot I_{1,n}(t)\\ \svdots \\ \dot I_{m,1}(t) \\ \svdots \\ \dot I_{m,n}(t) \end{smallpmatrix} =
\underbrace{\begin{bmatrix} \tilde A & & \\ & \ddots & \\ & & \tilde A\end{bmatrix}}_{=A}
\underbrace{\begin{smallpmatrix}I_{1,1}(t) \\ \svdots \\ I_{1,n}(t)\\ \svdots \\ I_{m,1}(t) \\ \svdots \\ I_{m,n}(t) \end{smallpmatrix}}_{=:\bar I(t)\in \R^{mn}}
+\underbrace{\begin{pmatrix}b_1(t) \\ \svdots \\ b_m(t)\end{pmatrix}}_{=:b(t)\in\R^{mn}}
\end{equation}
for $t\geq 0$, where for $\im$ we have $b_i(t)\in \R^n$ defined by
\begin{equation*}
b_i(t)
:=\bigg[0,\ldots,0, \sum\limits_{k=1}^m s_{i,k} \big(u_k(t)\!-\!u_k(t\!-\!\tus\!-\!\tau_u(t\!-\!\tus))\big)\bigg]^\top\!.
\end{equation*}
Defining $\zeta(t):= \eta(t\!-\!\tus)$ the last of equations~\eqref{eq:Sys} becomes
\[
    \dot \zeta(t) = h(t\!-\!\tus,\bar x_n(t\!-\!\tus),\zeta(t)).
\]
Then, by~\eqref{eq:Sys},~\eqref{eq:controller},~\eqref{eq:controller-I} and Step~1 there exists a function
\[H:\R_{\geq \tau_s}\times(-1,1)^{m(n-1)}\times \Omega^2\times \Omega_I^{mn}\times\R^q\to \R^{2mn+q}\] satisfying
\begin{equation}\label{eq:ODE_bar-z}
    \begin{smallpmatrix} \dot{\bar{z}}(t)\\ \dot{\bar{I}}(t)\\ \dot \zeta(t)\end{smallpmatrix}= H\big(t,\bar z(t), z_n(t\!-\!\tus\!-\!\tau_u(t\!-\!\tus)),\bar I(t),\zeta(t)\big),
\end{equation}
where $\bar{z}(t) =\big[z_{1,1}(t) \dots z_{m,1}(t) \dots z_{1,n}(t) \ \dots \ z_{m,n}(t)\big]^{\top}$.
Note that $H$ is piecewise continuous and bounded in~$t$ and locally Lipschitz in all other variables. We define
the operator
\begin{multline*}
    \scs:\cC([-\bar\tau_u,\infty),\R^{m})\to\mathcal{L}_\text{loc}^\infty(\R_{\geq \tau_s},\R^{n}),\\ \xi \mapsto \big(t\mapsto \xi(t\!-\!\tus\!-\!\tau_u(t\!-\!\tus))\big).
\end{multline*}
We can now rewrite~\eqref{eq:ODE_bar-z} as
\begin{equation}\label{eq:ODE_bar-z_new}
    \begin{smallpmatrix} \dot{\bar{z}}(t)\\ \dot{\bar{I}}(t)\\ \dot \zeta(t)\end{smallpmatrix}= \tilde H\big(t,\bar z(t), \bar I(t),\zeta(t), \tilde \scs(\bar z)(t)\big)
\end{equation}
for some 
\[\tilde H:\R_{\geq \tau_s}\times(-1,1)^{m(n-1)}\times \Omega\times\Omega_I^{mn}\times\R^{q+n}\to \R^{mn+q},\]
and where $\tilde \scs$ is defined such that $\tilde \scs(\bar z)=\scs(z_n)$. Furthermore, by assumption of the theorem, there exists a well-defined initial history function $z^\varphi\in \cC([-\bar\tau_u,\tau_s],\R^{mn})$ for $\bar z$ with $z^\varphi(t) \in (-1,1)^{m(n-1)}\times\Omega$ for all $t\in [-\bar\tau_u,\tau_s]$. Therefore, it follows from~\eqref{eq:controller} and the fact that $u_i(t) = 0$ for $t\in [-\tau_s-\bar\tau_u,0]$ and $\im$ that
\[
    u_i(t) = \begin{cases} \frac{-\sigma k_n \chi(\|z^\varphi_n(t)\|) z^\varphi_{i,n}(t)}{1-\norm{z^\varphi_n(t)}^2}, & t\in[0,\tau_s],\\
    0, & t\in [-\bar\tau_u,0]. \end{cases}
\]  
Hence, there exists a well-defined initial history function $I^\varphi\in \cC([0,\tau_s],\R^{mn})$ for $\bar I$, given by
\[
    I^\varphi(t) = \begin{cases} \int_0^t e^{A(t-s)} b^\varphi(s) \ds{s},& t\in [0,\tau_s],\\ 0,& t\in [-\bar\tau_u,0]\end{cases}
\]
with $b^\varphi(s) = (b_1^\varphi(s)^\top, \ldots, b_m^\varphi(s)^\top)^\top$ 
and $b_i^\varphi(s)=\left[0 \ \dots 0 \ \ \sum\limits_{k=1}^m s_{i,k} u_k(s)\right]^{\top}$.
Now, as a first lower bound on $\delta_I$, we assume that $\delta_I >  |I^\varphi_{i,j}(t)|$, for all $t\in[0,\tau_s]$, all $i=1,\ldots,m$ and all $j=1,\ldots,n$; this bound can be determined a priori.

The initial history for $\zeta$ is, according to~\eqref{eq:initial_history}, simply given by $\zeta^\varphi(\cdot) = \eta^\varphi(\cdot-\tus)\in \cC([-\bar\tau_u,\tau_s],\R^{q})$. Overall,~\eqref{eq:ODE_bar-z_new} is equipped with the initial history
\begin{equation}\label{eq:ODE_bar-z_new-hist}
    \big(\bar z,\bar I,\zeta)\vert_{[-\bar\tau_u,\tau_s]} = (z^\varphi, I^\varphi, \zeta^\varphi).
\end{equation}
The existence of a maximal solution $(\bar z,\bar I,\zeta): [-\bar\tau_u,\omega)\to \R^{2mn+q}$, $\omega\in( \tau_s,\infty]$, of~\eqref{eq:ODE_bar-z_new},~\eqref{eq:ODE_bar-z_new-hist} such that the closure of $\setdef{\big(t,\bar z(t), \bar I(t), \zeta(t)\big)}{t\in [\tau_s,\omega)}$
is not a compact subset of $\R_{\geq \tau_s}\times\big((-1,1)^{m(n-1)}\times\Omega\big)\times \Omega_I^{mn} \times\R^q$ then follows from~\cite[Thm.~B.1]{IlchRyan09}. 

\emph{Phase 3}: We show that each $z_{i,j}(t)$ evolves strictly within $(-1,1)$ and $z_n(t)$ evolve strictly within $\Omega$, $\im$,\linebreak $\jn-1$, for all $t\in[\tau_s,\omega)$. The proof follows a recursive procedure.

\textit{\textbf{Step 1 ($j=1$, $t\in [\tau_s,\omega)$):}} We begin by considering the positive functions $V_{i,1}(t):=\frac 1 2 z_{i,1}(t)^2$ for $t\in[ \tau_s,\omega)$ and $\im$. By Step~1 we have
\begin{equation}\label{eq:est-Vi1-1}
\begin{aligned}
&\dot V_{i,1}(t)
=\tfrac{z_{i,1}(t)}{\psi_{i,1}(t\!-\!\tus)}\Big(z_{i,2}(t)\psi_{i,2}(t\!-\!\tus)-\tfrac{k_{i,1}z_{i,1}(t)}{1-z_{i,1}(t)^2}\\
& 	\quad -\dot y_{d,i}(t\!-\!\tus)-\dot \psi_{i,1}(t\!-\!\tus)z_{i,1}(t) + f_{i,1}\big(t, \bar x_1(t\!-\!\tus)\big)\Big).
\end{aligned}
\end{equation}
Recall that $\bar x_1$ can be rewritten in terms of $z_{i,1}$, $\psi_{i,1}$, $y_{d,i}$ and $I_{i,1}$, as mentioned in Step~1, i.e., there exists a continuous function $F_1: \R^{4m}\to\R^{mn}$ such that, for all $t\ge \tau_s$, 
\begin{multline*}
    \bar x_1(t\!-\!\tus) = F_1\big(z_{1,1}(t),\ldots,z_{m,1}(t), I_{1,1}(t),\ldots,I_{m,1}(t),\\
    \psi_{1,1}(t),\ldots,\psi_{m,1}(t), y_{d,1}(t),\ldots,y_{d,m}(t)\big).
\end{multline*}
Since $z_{i,1}(t),z_{i,2}(t)\in(-1,1)$, $I_{i,1}(t)\in\Omega_I$, $y_{d,i}$, $\psi_{i,1}$ are bounded and $f_{i,1}$ are locally Lipschitz in~$\bar x_{i,1}$ and bounded in~$t$, we deduce that
\[
     |f_{i,1}\big(t, \bar x_1(t\!-\!\tus)\big)| \le \sup_{\xi\in\Omega_{\xi,1}}|f_{i,1}(t,F_1(\xi))|:=\bar{f}_{i,1},
\]
where
\begin{equation*}
    \Omega_{\xi,1}:=(-1,1)^m\times\Omega_I^m\times [-\bar \psi_1, \bar \psi_1]^m\times  [-\bar y_d, \bar y_d]^m
\end{equation*}
and $\bar \psi_1 = \max_{i=1,\ldots,m} \sup_{t\ge 0} |\psi_{i,1}(t)|$, $\bar y_d = \max_{i=1,\ldots,m} \sup_{t\ge 0} |y_{d,i}(t)|$. Moreover, $\dot \psi_{i,1}$ and $\psi_{i,2}$ are bounded by construction and $\dot y_{d,i}$ is bounded by assumption, we can utilize~\eqref{eq:est-Vi1-1} to obtain
\begin{equation}\label{eq:est-Vi1-2}
\dot V_{i,1}(t) \leq \tfrac{|z_{i,1}(t)|}{\psi_{i,1}(t\!-\!\tus)}\left(C_{i,1}-\tfrac{k_{i,1}|z_{i,1}(t)|}{1-z_{i,1}(t)^2}\right)
\end{equation}
for some $C_{i,1}>0$. Choose $\vareps_{i,1}\in(0,1)$ such that
\[
    \vareps_{i,1}>\max\left\{-\tfrac{k_{i,1}}{2C_{i,1}}+\sqrt{\left(\tfrac{k_{i,1}}{2C_{i,1}}\right)^2\!+\!1},\sup_{t\in[-\bar\tau_u, \tau_s]}|z_{i,1}(t)|\right\}.
\]
Assume that there exists $t_1\in( \tau_s,\omega)$ with $|z_{i,1}(t_1)|>\vareps_{i,1}$.
Since $|z_{i,1}(t)|\leq\vareps_{i,1}$ for all $t\in[-\bar\tau_u, \tau_s]$, $t_0:=\max \left\{t\in[ \tau_s,t_1)\,\big\vert\, |z_{i,1}(t)|=\vareps_{i,1}\right\}$
is well-defined. Then we have $|z_{i,1}(t)|>\vareps_{i,1}$ and $\frac{k_{i,1}|z_{i,1}(t)|}{1-z_{i,1}(t)^2}>\frac{k_{i,1}\vareps_{i,1}}{1-\vareps_{i,1}^2}$,
%
for all $t\in(t_0,t_1]$. Utilizing that for $a\in [0,1)$ we have
\[ C_{i,1}-\tfrac{k_{i,1} a}{1-a^2}  <  0 \
 \Leftrightarrow\  -\tfrac{k_{i,1}}{2C_{i,1}}+\sqrt{\left(\tfrac{k_{i,1}}{2C_{i,1}}\right)^2+1}  <  a,
\]
we find that, by construction of $\vareps_{i,1}$, $C_{i,1}-\frac{k_{i,1}|z_{i,1}(t)|}{1-z_{i,1}(t)^2} < 0$ for all $t\in (t_0,t_1]$.
Then it follows from~\eqref{eq:est-Vi1-2}  that $\dot V_{i,1}(t)<0$ for all $t\in (t_0,t_1]$.
Thus we arrive at the contradiction $\vareps_{i,1}=|z_{i,1}(t_0)|>|z_{i,1}(t_1)|>\vareps_{i,1}$.
Therefore, we have shown that
\begin{equation}\label{eq:zi1-epsi1}
\forall\, t\in[-\bar\tau_u,\omega):\ |z_{i,1}(t)|\leq \vareps_{i,1}.
\end{equation}
We once again recall Step~1 to infer that also $\dot z_{i,1}$ is bounded on $[-\bar\tau_u,\omega)$ for all $\im$. Then it follows from~\eqref{eq:controller} that $\dot a_{i,1}$ is bounded on $[-\bar\tau_u,\omega)$, $\im$.

\textit{\textbf{Step $j$ ($j=2,\ldots,n-1$, $t\in [\tau_s,\omega)$:}}
Again consider the positive functions $V_{i,j}(t):=\frac 1 2 z_{i,j}(t)^2$, for $t\in[ \tau_s,\omega)$ and $\im$. By Step~1, the only difference between the differential equations for $V_{i,1}$ and $V_{i,j}$, except from the changed indices, is the term $\dot y_{d,i}(t\!-\!\tus)$ appearing in $\dot V_{i,1}(t)$, which is replaced by $\dot a_{i,j-1}(t)$ in $\dot V_{i,j}(t)$. Furthermore, $\bar x_j(t\!-\!\tus)$ depends (continuously) on more of the previous variables (cf.\ Step~1), again all of which are bounded. Since $\dot a_{i,1}$ is bounded by the proof of the Case $j=1$, we may inductively show that there exists $\eps_{i,j}\in (0,1)$ such that
\begin{equation}\label{eq:zij-epsij}
\forall\, t\in[-\bar\tau_u,\omega):\ |z_{i,j}(t)|\leq \vareps_{i,j}.
\end{equation}
and $\dot z_{i,j}$ as well as $\dot a_{i,j}$ are bounded on $[-\bar\tau_u,\omega)$ for all $\im$ .

\textit{\textbf{Step $n$ ($j=n$, $t\in [\tau_s,\omega)$):}} We consider the positive function $V_n(t)	:=\frac 1 2 \norm{z_n(t)}^2 =\frac 1 2 \sum\limits_{i=1}^m z_{i,n}(t)^2$.
By~\eqref{eq:controller},~\eqref{eq:controller-I} and Step~1 we calculate the derivative of $V_n(t)$ as
\begin{equation}\label{eq:dotVn}
\begin{aligned}
&\dot V_n(t)
=z_n(t)^\top \begin{smallpmatrix}\dot z_{1,n}(t) \\ \svdots \\ \dot z_{m,n}(t)\end{smallpmatrix}\\
&=	\frac{z_n(t)^\top}{\psi_n(t\!-\!\tus)} \Bigg(\begin{smallpmatrix}f_1(t\!-\!\tus,\bar x_n(t\!-\!\tus),\eta(t\!-\!\tus))\\ \svdots \\ f_m(t\!-\!\tus,\bar x_n(t\!-\!\tus),\eta(t\!-\!\tus))\end{smallpmatrix} \\
& - \dot \psi_{n}(t\!-\!\tus) \begin{smallpmatrix}z_{1,n}(t) \\ \svdots \\z_{m,n}(t)\end{smallpmatrix}-\begin{smallpmatrix}\dot a_{1,n-1}(t) \\ \svdots \\ \dot a_{m,n-1}(t)\end{smallpmatrix} -\tilde I(t) +S\bar u(t)\\
&  +\big(G(t\!-\!\tus,\bar x_n(t\!-\!\tus),\eta(t\!-\!\tus)) 
 -S\big) \bar u(t\!-\!\tus\!-\!\tau_u(t\!-\!\tus)) \Bigg),
\end{aligned}
\end{equation}
where
\begin{align*}
\tilde I(t)&:=\begin{smallpmatrix}
\sum\limits_{k=1}^{n}\tbinom{n}{n+1-k}(-\alpha)^{n+1-k}I_{1,k}(t)\\
\svdots\\
\sum\limits_{k=1}^{n}\tbinom{n}{n+1-k}(-\alpha)^{n+1-k}I_{m,k}(t)
\end{smallpmatrix},
\end{align*}
and $\bar u(t):=[u_1(t) \ \dots \ u_m(t)]^{\top}$.
By recalling Assumption~\ref{Ass1} and~\eqref{eq:controller},~\eqref{eq:controller-I} we can derive the estimate 
\begin{equation}\label{eq:est-fi}
\small{
\begin{aligned}
&\vert f_{i,n}(t\!-\!\tus,\bar x_n(t\!-\!\tus),\eta(t\!-\!\tus))\vert \\
&\leq \vert d_i(t\!-\!\tus)\vert\big(\norm{\bar x_{n}(t\!-\!\tus)}+\norm{\eta(t\!-\!\tus)}+1\big)\\
&\leq \vert d_i(t\!-\!\tus)\vert\left(1+ \vphantom{\begin{pmatrix}a\\ \svdots   \\ a\end{pmatrix}}
\norm{\begin{smallpmatrix}\psi_{1,1}(t\!-\!\tus)z_{1,1}(t) \\ \svdots \\ \psi_{m,n}(t\!-\!\tus)z_{m,n}(t)\end{smallpmatrix}} \right.\\
& +\kappa(t\!-\!\tus,\norm{\eta(0)}) + \bar \gamma_1 \sup_{s\in[0,t\!-\!\tus]} \norm{\bar x_n(s)} + \bar \gamma_2 + c\\
&  + \norm{\begin{smallpmatrix}y_{d,1}(t\!-\!\tus) \\ \svdots \\ y_{d,m}(t\!-\!\tus) \\ a_{1,1}(t) \\ \svdots \\ a_{m,n-1}(t)\end{smallpmatrix}}+\left.\norm{\begin{smallpmatrix}I_{1,1}(t) \\ \svdots  \\ \sum\limits_{k=1}^j\tbinom{j-1}{j-k}(-\alpha)^{j-k}I_{i,k}(t)\\ \svdots \\ \sum\limits_{k=1}^n\tbinom{n-1}{n-k}(-\alpha)^{n-k}I_{m,k}(t)\end{smallpmatrix}}\right)\\
&\leq |d_i(t\!-\!\tus)|\Big(\tilde C_1 +\max\left\{\alpha^{n-1},1\right\}\tbinom{n}{\lfloor\frac{n}{2}\rfloor}n\big(\norm{\bar I(t)} \\
&\hspace{5.5cm}+ \bar \gamma_1 \sup_{s\in[\tus,t]} \norm{\bar I(s)}\big)\Big)\\
&\le C_1  (1+ \sup_{s\in[\tus,t]} \norm{\bar I(s)})
\end{aligned}
}
\end{equation}
where $\tilde C_1, C_1>0$ exist due to the boundedness of $\kappa(\cdot,\|\eta(0)\|)$, $\psi_{i,j}$, $z_{i,j}$, $y_{d,i}$, $d_i$, $\im$, $\jn$ and $a_{i,j}$, $\im$, $\jn-1$.

We like to emphasize at this point, that the above estimate is quite conservative and one might find a better one depending on the given system. As this estimate will essentially dictate the upper bound for $\tus+\tau_u(t\!-\!\tus)$ later on, a relaxation would allow for higher state measurement and control input delays.

To utilize~\eqref{eq:est-fi} in~\eqref{eq:dotVn} we need to find an estimate for~$\norm{\bar I(t)}$. We begin by defining $\tilde u_i(t):=\left[0 \ \dots \ 0 \ \sum\limits_{k=1}^m s_{i,k} u_{k}(t)\right]^{\top}\in\R^n$, $\im$, and $\tilde u(t):=[\tilde u_1(t) \ \dots \ \tilde u_m(t)]^{\top}\in\R^n$, $\im$.
for $t\in[\tau_s,\omega)$. Since $t\mapsto t\!-\!\tau_s\!-\!\tau_u(t\!-\!\tus)$ is strictly monotonically increasing on $\R_{\ge \tus}$ by Assumption~\ref{Ass2}, there exists a strictly monotonically increasing function $\Gamma:[-\tau_u(0),\infty)\to\R$ such that $\Gamma( t\!-\!\tau_s\!-\!\tau_u(t\!-\!\tus)) = t$ for all $t\ge \tus$. 
\blue{Furthermore, we have that $f(\Gamma(r)) = r$ for all $r\ge -\tau_u(0)$ and hence $\Gamma(r) \le r$ for all $r\ge -\tau_u(0)$.}
By~\eqref{eq:ODE_Iij} and the fact that $u_i(t) = 0$ for $t\in [-\tau_s-\bar\tau_u,0]$ and $\im$,  \textit{variation of constants} leads to
\small{
\begin{align*}
\bar I(t)
&= \int_{t\!-\!\tau_s\!-\!\tau_u(t\!-\!\tus)}^t e^{A(t-s)} \tilde u(s) \intd s \\
&\quad - \int_{0}^{t\!-\!\tau_s\!-\!\tau_u(t\!-\!\tus)} \left(I - \tfrac{e^{A(r-\Gamma(r))}}{1-\dot\tau_u(\Gamma(r))}\right) e^{A(t-r)}\tilde u(r) \intd r.
\end{align*}
}
We can now estimate $\norm{\bar I(t)}$ by employing Assumption 3 and the estimates~\eqref{eq:est-eAt}, $\mu \le \|A\|$, $|r-\Gamma(r)| \le \tau_s + \bar\tau_u$ for $r\ge 0$, \blue{$1-e^{-y}\le y$ for $y\ge 0$, and $\|I-e^{Az}\| \le \int_0^z \|A e^{As}\| \intd s \le \int_0^z \|A\| M e^{-\mu s}\intd s
= \|A\|\frac{M}{\mu} (1-e^{-\mu z}) \le \|A\| M z$
for all $z\ge 0$,} so that \bluen{
\begin{equation}\label{eq:est-barI}
\small{
\norm{\bar I(t)}
\le \|S\|\, C_\tau \sup_{s\in[0,t]}\|\bar u(s)\|
}
\end{equation}
with $C_\tau$ as in Assumption~\ref{Ass3}.}
Now observe that $\bar u(t) = -\sigma k_n \chi(\|z_n(t)\|) \frac{z_n(t) }{1-\norm{z_n(t)}^2}$.
Using this and~\eqref{eq:est-fi} in combination with the boundedness of $\dot \psi_n$, $z_{i,n}$ and $\dot a_{i,n-1}$, $\im$ as well as $\Theta:= \sup_{(t,x,\eta)\in\R_{\ge 0}\times\R^{nm+q}} \norm{G(t,x,\eta)-S} < \infty$,
which exists by Assumption 3, in~\eqref{eq:dotVn} we arrive at the inequality
\begin{align*}
&\dot V_n(t)
\leq \blue{\tfrac{\norm{z_n(t)}}{\psi_n(t-\tau_s)} } \bigg( \blue{k_n} C_3 + C_4 \sup_{s\in[\tus,t]} \norm{\bar I(s)} \\
& +  \Theta  \norm{\bar u(t\!-\!\tus\!-\!\tau_u(t\!-\!\tus))}- \blue{k_n \chi(\|z_n(t)\|)} \tfrac{s^*\norm{z_n(t)}}{1-\norm{z_n(t)}^2}\bigg),
\end{align*}
for some $C_2,C_3,C_4>0$, where $s^*$ is defined in Sec.~\ref{Sec:ContrDes} and
\begin{multline*}
    C_4 := (1+\bar \gamma_1) \max\left\{\alpha^{n-1},1\right\}\tbinom{n}{\lfloor\frac{n}{2}\rfloor} n\cdot \max_{\im}\|d_i\|_\infty \\
    +\max\left\{\alpha^{n},1\right\}\tbinom{n}{\lfloor\frac{n}{2}\rfloor}.
\end{multline*}
Next we employ~\eqref{eq:est-barI} to further estimate 
\begin{equation}\label{eq:est-dotVn}
\begin{aligned}
&\dot V_n(t)
\leq \blue{\tfrac{\norm{z_n(t)}}{\psi_n(t-\tau_s)} }\bigg(\blue{k_n} C_3 + C_5 \sup_{s\in[0,t]}\tfrac{\blue{k_n (1-\delta)}\norm{z_n(s)}}{1-\norm{z_n(s)}^2} \\
& + \Theta \norm{\bar u(t\!-\!\tus\!-\!\tau_u(t\!-\!\tus))} -\blue{k_n \chi(\|z_n(t)\|)} \tfrac{s^*\norm{z_n(t)}}{1-\norm{z_n(t)}^2}\bigg),
\end{aligned} 
\end{equation}
for \bluen{$C_5 := C_4\,\|S\|\,C_\tau$,} according to~\eqref{eq:est-barI}. We note that $\norm{\bar u(t\!-\!\tus\!-\!\tau_u(t\!-\!\tus))} \\
    = \frac{k_n \chi(\|z_n(t\!-\!\tus\!-\!\tau_u(t\!-\!\tus))\|)\norm{z_n(t\!-\!\tus\!-\!\tau_u(t\!-\!\tus))}}{1-\norm{z_n(t\!-\!\tus\!-\!\tau_u(t\!-\!\tus))}^2}$ for $t\ge \tus +\tau_u(t\!-\!\tus)$ and $\norm{\bar u(t\!-\!\tus\!-\!\tau_u(t\!-\!\tus))} =0$ otherwise. Choose $\varepsilon_n\in(0,1)$ with $\varepsilon_n > \max\bigg\{-\tfrac{C_6}{2C_3}+\sqrt{\left(\tfrac{C_6}{2C_3}\right)^2+1},\sup_{t\in[-\bar\tau_u, \tau_s]}\norm{z_n(t)}, \tfrac{1+\delta}{2} \bigg\}$,
where $\delta\in [0,1)$ is as in~\eqref{eq:act_fct} and $C_6 := \blue{\frac{1-\delta}{2}s^*- (1-\delta) C_5- \Theta}  > 0$,
by Assumption 3. Assume that there exist $\zeta\in (\varepsilon_n,1)$ and $t\in(\tau_s,\omega)$ with $\norm{z_n(t)}>\zeta$. Define $t_0:= \inf \setdef{t\in[\tau_s,\omega)}{\|z_n(t)\| > \zeta}$
which is well-defined since $\norm{z_n(t)}\le \varepsilon_n$ for all $t\in[-\bar\tau_u, \tau_s]$. By continuity there exists $t_1\in (t_0,\omega)$ such that $\|z_n(t_1)\| > \zeta$ and $\|z_n(t)\| \ge \vareps_n > \delta$ for all $t\in [t_0,t_1]$. Therefore, we find that
\begin{multline}\label{eq:chi-delta}
   \forall\, t\in[t_0,t_1]:\ \chi(\|z_n(t)\|) \ge \|z_n(t)\|-\delta \\
   \ge \eps_n-\delta > \tfrac{1+\delta}{2} -\delta = \tfrac{1-\delta}{2}.
\end{multline}
Furthermore, we observe that 
\begin{equation}\label{eq:norm-zn-t0}
 \forall\, t\in[-\bar\tau_u, t_0]:\ \norm{z_n(t)} \le \zeta = \norm{z_n(t_0)}.
\end{equation}
Next we distinguish two cases.

\textit{Case I:} There exists $t_2\in (t_0,t_1]$ such that $\dot V_n(t) \le 0$ for all $t\in [t_0,t_2)$. Clearly, $t_2$ can be chosen such that $\|z_n(t_2)\|>\zeta$, otherwise this contradicts the definition of~$t_0$. Then we directly obtain the contradiction $\zeta = \|z_n(t_0)\| \ge \|z_n(t_2)\| > \zeta$.

\textit{Case II:} Assuming the opposite of Case~I leads to existence of a sequence $(t_k)\subset (t_0,t_1)$ with $t_k\to t_0$ for $k\to \infty$ such that $\dot V_n(t_k) > 0$ for all $k\in\N$. For $k\in\N$, define
$
    \tau_k := \max\{t\in [t_0,t_k) \mid \dot V_n(t) = 0\},
$
then we have $\dot V_n(t) \ge 0$ for all $t\in [\tau_k,t_k]$ and all $k\in\N$. Choose $\rho>0$ small enough so that
\begin{equation}\label{eq:delta}
 \vareps_n > -\tfrac{C_6}{2(C_3 + \rho\blue{(1-\delta)} C_5)}+\sqrt{\left(\tfrac{C_6}{2(C_3 + \rho\blue{(1-\delta)} C_5)}\right)^2+1},
\end{equation}
which is always possible. Then, by continuity, there exists $k\in\N$ sufficiently large such that
$
     \sup_{s\in[t_0,\tau_k]}\tfrac{\norm{z_n(s)}}{1-\norm{z_n(s)}^2} \le \tfrac{\norm{z_n(\tau_k)}}{1-\norm{z_n(\tau_k)}^2} + \rho.
$
Together with~\eqref{eq:norm-zn-t0} and monotonicity of $V_n$ on~$[\tau_k,t_k]$ this implies that
\begin{equation}\label{eq:zn-mono}
    \forall\, t\in[\tau_k,t_k]:\ \sup_{s\in[0,t]}\tfrac{\norm{z_n(s)}}{1-\norm{z_n(s)}^2} \le \tfrac{\norm{z_n(t)}}{1-\norm{z_n(t)}^2} + \rho.
\end{equation}
Incorporating~\eqref{eq:chi-delta} and~\eqref{eq:zn-mono} in~\eqref{eq:est-dotVn} leads to
\begin{align*}
\dot V_n(t)
&<\blue{\tfrac{k_n \varepsilon_n}{\psi_n(t-\tau_s)}}\left( C_3  + \rho \blue{(1-\delta)} C_5 - C_6\tfrac{\varepsilon_n}{1-\varepsilon_n^2}\right)<0,
\end{align*}
where for the last step we used~\eqref{eq:delta} and the fact that for any $a\in [0,1)$
\[\begin{array}{crcl}
 & (C_3 + \rho \blue{(1\!-\!\delta)} C_5)-C_6\frac{a}{1-a^2} & < & 0\\
 \Leftrightarrow & \frac{-C_6}{2(C_3 + \rho \blue{(1\!-\!\delta)} C_5)}+\sqrt{\left(\frac{C_6}{2(C_3 + \rho \blue{(1\!-\!\delta)} C_5)}\right)^2+1} & < & a.
\end{array}\]
Now $\dot V_n(t)<0$ for all $t\in[\tau_k,t_k]$ directly contradicts the choice of the interval $[\tau_k,t_k]$. Overall, we have shown that
\begin{equation}\label{eq:zn-epsn}
    \forall\, t\in[-\bar\tau_u,\omega):\ \|z_n(t)\| \le \vareps_n.
\end{equation}

\emph{Phase 4}: Employing \eqref{eq:zn-epsn}, we further conclude the existence of constants $\bar{u}_i>0$, $i=1,\hdots,n$, which are independent of $\delta_I$ (since $C_3$ and $C_6$ which define $\eps_n$ are independent of $\delta_I$), such that $|u_i(t)|\leq\bar{u}_i$, $i=1,\hdots,m$, for all $t\in[\tau_s,\omega)$. Therefore, invoking the definition of $I_{i,j}$, we conclude that we may choose $\delta_I$ sufficiently large such that for all $i=1,\hdots,m$, $j=1,\hdots,n$, and $t\in[\tau_s,\omega)$, 
\begin{equation*}
    |I_{i,j}(t)|\leq \sum_{\ell=1}^{m}\frac{2s_{i,\ell}\bar{u}_i}{\alpha^{n-j}}+\sum_{w=j+1}^{n}\frac{|I_{i,j}^\varphi(\tau_s)|}{\alpha^{w-j}}+|I_{i,j}^\varphi(\tau_s)|<\delta_I,
\end{equation*}
guaranteeing that $I_{i,j}$ evolves strictly within a compact subset of $\Omega_I$ for all $t\in[\tau_s,\omega)$. Furthermore, it follows from~\eqref{eq:zi1-epsi1},~\eqref{eq:zij-epsij}, and~\eqref{eq:zn-epsn} that $\bar z$ evolves in a compact subset of $(-1,1)^{m(n-1)}\times \Omega$. Therefore, Step~2 implies that $\omega=\infty$.
Furthermore, we have shown that all closed-loop signals remain bounded. \bluen{Finally, (iv) follows from (iii), $|I_{i,1}(t)|\le\|\bar I(t)\|$
and~\eqref{eq:est-barI}; note that \eqref{eq:est-barI} also holds for
$t\in[0,\tau_s]$, since
$\|I^\varphi(t)\| \le M\tau_s\|S\|\sup_{s\in[0,t]}\|\bar u(s)\|
 \le \|S\|C_\tau\sup_{s\in[0,t]}\|\bar u(s)\|$.} This finishes the proof of the theorem.\hfill $\square$

\bibliographystyle{IEEEtran}
\bibliography{references}

@STRING{ ijc        = {Int. J. Control} }

@STRING{ mm         = {Manuscripta Mathematica} }

@STRING{ mcss       = {Math. Control Signals Syst.} }

@STRING{ siamjco    = {{SIAM} J. Control Optim.} }

@STRING{ scl        = {Syst. Control Lett.} }

@article{BergDenn24,
  author = {Berger, Thomas and Dennst\"{a}dt, Dario and Lanza, Lukas and Worthmann, Karl},
  title = {Robust {F}unnel {M}odel {P}redictive {C}ontrol for output tracking with prescribed performance},
  journal=siamjco,
  year = {2024},
  volume = 62,
  number = 4,
  pages = {2071--2097},
}

@article{LanzDenn24,
  title={Sampled-data funnel control and its use for safe continual learning},
  author={Lanza, Lukas and Dennst{\"a}dt, Dario and Worthmann, Karl and Schmitz, Philipp and \c{S}en, Gökçen Devlet and  Trenn, Stephan and Schaller, Manuel},
  journal=scl,
  year={2024},
  volume = 192, 
  pages = {Article 105892},
}

@ARTICLE{BergIlch21,
   AUTHOR    = {Berger, Thomas and Ilchmann, Achim and Ryan, Eugene P.},
   YEAR      = 2021,
   TITLE     = {Funnel control of nonlinear systems},
   JOURNAL   = mcss,
   volume    = 33,
   pages     = {151--194},
}

@ARTICLE{IlchRyan02b,
   AUTHOR    = {Ilchmann, Achim and Ryan, Eugene P. and Sangwin, Christopher J.},
   YEAR      = 2002,
   TITLE     = {Tracking with prescribed transient behaviour},
   JOURNAL   = {ESAIM: Control, Optimisation and Calculus of Variations},
   Volume    = 7,
   Pages     = {471--493},
   mstnote   = {MM, added 28.03.2008},
}

@ARTICLE{SeifBlaj13,
   AUTHOR    = {Seifried, Robert and Blajer, Wojciech},
   YEAR      = 2013,
   TITLE     = {Analysis of Servo-Constraint Problems for Underactuated Multibody Systems},
   JOURNAL   = {Mech. Sci.},
   Volume    = 4,
   Pages     = {113--129},
}

@ARTICLE{IlchRyan09,
   AUTHOR    = {Ilchmann, Achim and Ryan, Eugene P.},
   YEAR      = 2009,
   TITLE     = {Performance funnels and tracking control},
   JOURNAL   = ijc,
   Volume    = 82,
   Number    = 10,
   Pages     = {1828--1840},
}

@ARTICLE{BikaRovi23,
  author = {Bikas, L. N. and Rovithakis, G. A.},
  title = {Prescribed performance tracking of uncertain {MIMO} nonlinear systems in the presence of delays},
  journal = {{IEEE} Trans. Autom. Control},
  year = {2023},
  volume = 68,
  number = 1,
  pages = {96--107},
}

@Article{BergLeReis2018,
  author   = {Thomas Berger and Huy Hoàng Lê and Timo Reis},
  title    = {Funnel control for nonlinear systems with known strict relative degree},
  journal  = {Automatica},
  year     = {2018},
  volume   = {87},
  pages    = {345-357},
}

@article{RoviSurvey,
    author = {Rovithakis, George A.},
    title = {Prescribed Performance Adaptive Control of Uncertain Nonlinear Systems: State-of-the-art and Open Issues},
    journal = {PAMM},
    volume = {18},
    number = {1},
    pages = {e201800134},
    year = {2018}
}

@ARTICLE{Bech2008,
  author={Bechlioulis, Charalampos P. and Rovithakis, George A.},
  journal = {{IEEE} Trans. Autom. Control},
  title={Robust Adaptive Control of Feedback Linearizable MIMO Nonlinear Systems With Prescribed Performance}, 
  year={2008},
  volume={53},
  number={9},
  pages={2090-2099},
}

@ARTICLE{Karafyllis2012,
  author={Karafyllis, Iasson and Krstic, Miroslav},
  journal = {{IEEE} Trans. Autom. Control},
  title={Nonlinear Stabilization Under Sampled and Delayed Measurements, and With Inputs Subject to Delay and Zero-Order Hold}, 
  year={2012},
  volume={57},
  number={5},
  pages={1141-1154}
}

@article{Selivanov2016,
    title = {Predictor-based networked control under uncertain transmission delays},
    journal = {Automatica},
    volume = {70},
    pages = {101-108},
    year = {2016},
    author = {Anton Selivanov and Emilia Fridman},
}

@article{Zhou2017,
title = {Stabilization of linear systems with both input and state delays by observer–predictors},
journal = {Automatica},
volume = {83},
pages = {368-377},
year = {2017},
author = {Bin Zhou and Qingsong Liu and Frédéric Mazenc},
}

@article{Le2018,
title = {Prediction-based control of LTI systems with input and output time-varying delays},
journal = {Systems and Control Letters},
volume = {112},
pages = {24-30},
year = {2018},
author = {V. Léchappé and E. Moulay and F. Plestan},
}

@ARTICLE{Weston2019,
  author={Weston, Jerome and Malisoff, Michael},
  journal = {{IEEE} Trans. Autom. Control},
  title={Sequential Predictors Under Time-Varying Feedback and Measurement Delays and Sampling}, 
  year={2019},
  volume={64},
  number={7},
  pages={2991-2996},
}

@ARTICLE{Battilotti2020,
  author={Battilotti, Stefano},
  journal = {{IEEE} Trans. Autom. Control},
  title={Continuous-Time and Sampled-Data Stabilizers for Nonlinear Systems With Input and Measurement Delays}, 
  year={2020},
  volume={65},
  number={4},
  pages={1568-1583},
}

@article{Nozari2020,
title = {Event-triggered stabilization of nonlinear systems with time-varying sensing and actuation delay},
journal = {Automatica},
volume = {113},
pages = {108754},
year = {2020},
author = {Erfan Nozari and Pavankumar Tallapragada and Jorge Cortés},
}

@ARTICLE{Zhao2021,
  author={Zhao, Congran and Lin, Wei},
  journal = {{IEEE} Trans. Autom. Control},
  title={Global Stabilization by Memoryless Feedback for Nonlinear Systems With a Limited Input Delay and Large State Delays}, 
  year={2021},
  volume={66},
  number={8},
  pages={3702-3709},
}

@ARTICLE{Sun2021,
  author={Sun, Jiwei and Lin, Wei},
  journal = {{IEEE} Trans. Autom. Control},
  title={A Dynamic Gain-Based Saturation Control Strategy for Feedforward Systems With Long Delays in State and Input}, 
  year={2021},
  volume={66},
  number={9},
  pages={4357-4364},
}

@article{Yu2025,
title = {Universal output feedback control of a class of uncertain nonlinear systems with unknown delays in state, input and output},
journal = {Automatica},
volume = {174},
pages = {112135},
year = {2025},
author = {Xin Yu and Wei Lin},
}

@ARTICLE{Liberzon2013,
  author={Liberzon, Daniel and Trenn, Stephan},
  journal = {{IEEE} Trans. Autom. Control},
  title={The Bang-Bang Funnel Controller for Uncertain Nonlinear Systems With Arbitrary Relative Degree}, 
  year={2013},
  volume={58},
  number={12},
  pages={3126-3141},
}

@ARTICLE{Heemels2010,
  author={Heemels, W. P. Maurice H. and Teel, Andrew R. and van de Wouw, Nathan and Nešić, Dragan},
  journal = {{IEEE} Trans. Autom. Control},
  title={Networked Control Systems With Communication Constraints: Tradeoffs Between Transmission Intervals, Delays and Performance}, 
  year={2010},
  volume={55},
  number={8},
  pages={1781-1796}
}

\end{document}